# An isoperimetric inequality on the $\ell_p$ balls

## Sasha Sodin[*]

School of Mathematics, Raymond and Beverly Sackler Faculty of Exact Sciences, Tel Aviv University, Tel Aviv, 69978, Israel. E-mail: sodinale@tau.ac.il

**Abstract.** The normalised volume measure on the $\ell_p^n$ unit ball ($1 \leq p \leq 2$) satisfies the following isoperimetric inequality: the boundary measure of a set of measure $a$ is at least $cn^{1/p}\tilde{a}\log^{1-1/p}(1/\tilde{a})$, where $\tilde{a} = \min(a, 1-a)$.

**Résumé.** Nous prouvons une inégalité isopérimétrique pour la mesure uniforme $V_{p,n}$ sur la boule unité de $\ell_p^n$ ($1 \leq p \leq 2$). Si $V_{p,n}(A) = a$, alors $V_{p,n}^+(A) \geq cn^{1/p}\tilde{a}\log^{1-1/p}1/\tilde{a}$, où $V_{p,n}^+$ est la mesure de surface associée à $V_{p,n}$, $\tilde{a} = \min(a, 1-a)$ et $c$ est une constante absolue.

En particulier, les boules unités de $\ell_p^n$ vérifient la conjecture de Kannan–Lovász–Simonovits (*Discrete Comput. Geom.* **13** (1995)) sur la constante de Cheeger d'un corps convexe isotrope.

La démonstration s'appuie sur les inégalités isopérimétriques de Bobkov (*Ann. Probab.* **27** (1999)) et de Barthe–Cattiaux–Roberto (*Rev. Math. Iberoamericana* **22** (2006)), et utilise la représentation de $V_{p,n}$ établie par Barthe–Guédon–Mendelson–Naor (*Ann. Probab.* **33** (2005)) ainsi qu'un argument de découpage.

## 1. Introduction

We study the isoperimetric properties of the normalised volume measure

$$V_{p,n} = \text{Vol}\,|_{B_p^n} / \text{Vol}(B_p^n)$$

on the $\ell_p^n$ unit ball

$$B_p^n = \{x = (x_1, \ldots, x_n) \in \mathbb{R}^n \mid \|x\|_p^p = |x_1|^p + \cdots + |x_n|^p \leq 1\}, \quad 1 \leq p \leq 2.$$

Recall that the lower Minkowski content $\mu^+$ associated to a measure $\mu$ is defined as

$$\mu^+(A) = \liminf_{\varepsilon \to +0} \frac{\mu\{x \mid \text{dist}(x, A) \leq \varepsilon\} - \mu(A)}{\varepsilon}$$

---

[*]Supported in part by the European network PHD, FP6 Marie Curie Actions, RTN, Contract MCRN–511953.

for measurable sets $A$; we are interested in the behaviour of the isoperimetric function

$$\mathcal{I}_\mu(a) = \inf_{a \leq \mu(A) < 1/2} \mu^+(A) \tag{1}$$

for $\mu = V_{p,n}$.

**Theorem 1.** *There exists a universal constant $c > 0$ such that for $1 \leq p \leq 2$, $0 < a < 1/2$*

$$\mathcal{I}_{V_{p,n}}(a) \geq c n^{1/p} a \log^{1-1/p} \frac{1}{a}. \tag{2}$$

This extends the previously known case $p = 2$ (the Euclidean ball), for which Burago and Maz'ya [16] (see also [1, 13]) have found the complete solution to the isoperimetric problem (the extremal sets are half-balls and intersections of the ball with another one with orthogonal boundary). For $p = 1$, the theorem answers a question of Sergey Bobkov (cf. [9]).

The inequality (2) (and hence also (3) and (5) below) is sharp, up to a constant factor, unless $a$ is exponentially small (in the dimension $n$). In fact, the left- and right-hand sides of (5) are of the same order when

$$A = A_{\xi,t} = \{x \in \mathbb{R}^n \mid \langle x, \xi \rangle \geq t\}$$

is a coordinate half-space. In the last section of this note we explain how to obtain a sharp bound for smaller values of $a$ as well.

***Remark 1.*** *Not surprisingly, the complementary bound*

$$\mathcal{J}_{V_{p,n}}(a) \geq c n^{1/p} a \log^{1-1/p} \frac{1}{a} \tag{3}$$

*for*

$$\mathcal{J}_\mu(a) = \inf_{1/2 < \mu(A) < 1-a} \mu^+(A) \tag{4}$$

*is also true; in fact, $\mathcal{J}_\mu \equiv \mathcal{I}_\mu$ according to Proposition A (note however the slight asymmetry between the definitions of $\mathcal{I}_\mu$ and $\mathcal{J}_\mu$). Then, (2) and (3), combined with a trivial approximation argument for the case $V_{p,n}(A) = 1/2$, yield*

$$V_{p,n}^+(A) \geq c n^{1/p} (V_{p,n}(A) \wedge (1 - V_{p,n}(A))) \log^{1-1/p} \frac{1}{V_{p,n}(A) \wedge (1 - V_{p,n}(A))} \tag{5}$$

*(valid for any $A \subset \mathbb{R}^n$, with the common convention $0 \times \log 0^{-1} = 0$, $\vee = \max$ and $\wedge = \min$).*

***Remark 2.*** *The following definition of isoperimetric function is more common than (1), (4):*

$$\mathfrak{Is}_\mu(a) = \inf_{\mu(A) = a} \mu^+(A). \tag{6}$$

*Clearly, Theorem 1 is equivalent to (5) and to the inequality*

$$\mathfrak{Is}_{V_{p,n}}(a) \geq c n^{1/p} (a \wedge (1-a)) \log^{1-1/p} \frac{1}{(a \wedge (1-a))}.$$

Schechtman and Zinn [24] proved (in particular) the following estimate on the tails of the Euclidean norm with respect to $V_{p,n}$.

**Theorem (Schechtman–Zinn).** *There exist universal constants $c$ and $t_0$ such that the inequality*

$$V_{p,n}\{\|x\|_2 \geq t\} \leq \exp(-cnt^p) \tag{7}$$

*holds for $t \geq n^{-(2-p)/(2p)} t_0$, $1 \leq p \leq 2$.*

In the subsequent work [25], they showed that this inequality is a special case of a general concentration inequality. Recall that a median $\operatorname{Med} F$ of a real function $F$ on a probability space $(X, \mu)$ is (non-uniquely) defined by the inequalities

$$\mu\{F \geq \operatorname{Med} F\} \geq \frac{1}{2}, \qquad \mu\{F \leq \operatorname{Med} F\} \geq \frac{1}{2}.$$

**Theorem (Schechtman–Zinn).** *For any 1-Lipschitz function $F : B_p^n \to \mathbb{R}$ (meaning that $|F(x) - F(y)| \leq \|x - y\|_2$ for $x, y \in B_p^n$),*

$$V_{p,n}\{F > \operatorname{Med} F + t\} \leq C \exp(-c_1 n t^p), \quad 0 < t < +\infty. \tag{8}$$

Let us sketch the standard argument that recovers (8) (with $C = 1/2$ and $c_1 = c^p/p^p$) from (5). Let

$$\phi(h) = V_{p,n}\{F > \operatorname{Med} F + h\}, \quad h \geq 0.$$

Then (5) (applied to $A = \{F > \operatorname{Med} F + h\}$) implies

$$\phi'(h) \leq -cn^{1/p} \phi(h) \log^{1-1/p} \frac{1}{\phi(h)}$$

(where strictly speaking $\phi'$ stands for the left upper derivative). Therefore the inverse function $\psi = \phi^{-1}$ satisfies

$$\psi\left(\frac{1}{2}\right) = 0, \qquad \psi'(u) \geq -\left[cn^{1/p} u \log^{1-1/p} \frac{1}{u}\right]^{-1};$$

hence

$$\psi(u) = -\int_u^{1/2} \psi'(u_1) \, du_1 \leq \frac{\log^{1/p}(1/u) - \log^{1/p} 2}{cn^{1/p}/p} \leq \left[\frac{\log(1/2u)}{c_1 n}\right]^{1/p}$$

and $\phi(h) \leq \exp(-c_1 n h^p)/2$ as in (8).

The proof of Theorem 1 splits into two cases. In Section 2 we apply Bobkov's isoperimetric inequality [6] to deal with $a$ as small as $\exp(-Cn^{p/2})$.

For larger $a$, the representation of $V_{p,n}$ put forward by Barthe, Guédon, Mendelson and Naor [4] allows us to reduce the question to an analogous one for a certain product measure. We deal with this case in Section 3, making use of the Barthe–Cattiaux–Roberto isoperimetric inequality [3].

We devote the last section to remarks and comments.

## 2. Small sets

This section is based on the following isoperimetric inequality, due to Sergey Bobkov [6]. Recall that a probability measure $\mu$ on $\mathbb{R}^n$ is called log-concave if for any $A, B \subset \mathbb{R}^n$ and for any $0 < t < 1$

$$\mu((1-t)A + tB) \geq \mu(A)^{1-t} \mu(B)^t.$$

**Theorem (Bobkov).** *Let $\mu$ be a log-concave probability measure on $\mathbb{R}^n$. Then, for any $A \subset \mathbb{R}^n$ and any $r > 0$,*

$$\mu^+(A) \geq \frac{1}{2r}\left\{\mu(A)\log\frac{1}{\mu(A)} + (1-\mu(A))\log\frac{1}{1-\mu(A)} + \log\mu\{\|x\|_2 \leq r\}\right\}. \tag{9}$$

Apply (9) to $\mu = V_{p,n}$, which is log-concave according to the Brunn–Minkowski inequality. If $r > n^{-(2-p)/(2p)}t_0$, the Schechtman–Zinn theorem (7) yields

$$V_{p,n}\{\|x\|_2 \geq r\} \leq \exp(-cnr^p)$$

and

$$\log V_{p,n}\{\|x\|_2 \leq r\} \geq \log(1 - \exp(-cnr^p)) \geq -C\exp(-cnr^p)$$

(here and further $C$, $c$, $c_1$, $c_2$, $C'$, etc. denote universal constants that may change their meaning from line to line, unless explicitly stated).

On the other hand, for $V_{p,n}(A) < c'$,

$$(1 - V_{p,n}(A))\log\frac{1}{1 - V_{p,n}(A)} \geq c''V_{p,n}(A).$$

Hence for

$$r = \frac{1}{c^{1/p}n^{1/p}}\log^{1/p}\frac{C}{c''V_{p,n}(A)}$$

the sum of the last two terms in the right-hand side of (9) is not negative. We conclude:

**Proposition 1.** *There exist two universal constants $c, C > 0$ such that*

$$V_{p,n}^+(A) \geq cn^{1/p}V_{p,n}(A)\log^{1-1/p}\frac{1}{V_{p,n}(A)}$$

*for sets $A \subset \mathbb{R}^n$ such that $V_{p,n}(A) < \exp(-Cn^{p/2})$.*

### 3. Big sets

To complete the proof of Theorem 1, we need

**Proposition 2.** *There exists a constant $c' > 0$ such that*

$$V_{p,n}^+(A) \geq c'n^{1/p}V_{p,n}(A)\log^{1-1/p}\frac{1}{V_{p,n}(A)}$$

*for sets $A \subset \mathbb{R}^n$ such that $\exp(-Cn^{p/2}) \leq V_{p,n}(A) < 1/2$ (where $C$ is the same as in Proposition 1).*

Consider the product measure

$$\mu_{p,n} = \mu_p^{\otimes n} \otimes \nu_p,$$

where

$$d\mu_p(t) = \frac{\exp(-|t|^p)}{2\Gamma(1 + 1/p)}dt,$$

$$d\nu_p(t) = pt^{p-1}\exp(-t^p)\mathbf{1}_{[0,+\infty)}(t)dt;$$

define $T(z) = x/\|z\|_p$ for $z = (x,y) \in \mathbb{R}^{n+1} = \mathbb{R}^n \times \mathbb{R}$. Barthe, Guédon, Mendelson and Naor proved [4] that the map $T: \mathbb{R}^{n+1} \to \mathbb{R}^n$ pushes $\mu_{p,n}$ forward to $V_{p,n}$ (that is, $V_{p,n}(A) = \mu_{p,n}(T^{-1}A)$ for $A \subset \mathbb{R}^n$).

**Key fact (Barthe–Cattiaux–Roberto [3]).** *The measure $\mu = \mu_{p,n}$ satisfies the isoperimetric inequality*

$$\mathcal{I}_\mu(a) \geq ca \log^{1-1/p} \frac{1}{a}, \quad 0 < a < \frac{1}{2}. \tag{10}$$

The notation here differs slightly from that in [3]; therefore we provide a short explanation, reproducing the argument in the proof of [3], Theorem 21.

*Explanation.* First, the measures $\mu_p$ and $\nu_p$ are log-concave. Therefore we can use the following proposition.

**Proposition (Bobkov [10]).** *Let $\mu$ be a log-concave measure on $\mathbb{R}$; denote $F_\mu(x) = \mu(-\infty, x]$. Then*

$$\mathfrak{Is}_\mu(a) = \min(F'_\mu(F_\mu^{-1}(a)), F'_\mu(F_\mu^{-1}(1-a))), \quad 0 < a < 1$$

*(in other words, the extremal sets are half-lines).*

Now a simple computation shows that

$$c_1^{-1} \mathfrak{Is}_{\mu_p}(a) \geq (a \wedge (1-a)) \log^{1-1/p} \frac{1}{(a \wedge (1-a))} \geq c_1 \mathfrak{Is}_{\mu_p}(a);$$

$$\mathfrak{Is}_{\nu_p}(a) \geq c_2(a \wedge (1-a)) \log^{1-1/p} \frac{1}{(a \wedge (1-a))} \tag{11}$$

and thereby $\mathfrak{Is}_{\mu_p}, \mathfrak{Is}_{\nu_p} \geq c_3 \mathfrak{Is}_{\mu_p}$.

The measure $\mu_p$ is log-concave; therefore, according to the comparison theorem due to Barthe [2], Theorem 10,

$$\mathfrak{Is}_{\mu_{p,n}} = \mathfrak{Is}_{\mu_p^n \otimes \nu_p} \geq c_3 \mathfrak{Is}_{\mu_p^{n+1}}.$$

Barthe, Cattiaux and Roberto proved that the products of the measure $\mu_p$ satisfy a dimension-free isoperimetric inequality. More precisely, according to inequality (4) in the Introduction to [3],

$$\mathfrak{Is}_{\mu_p^{n+1}}(a) \geq c_4(a \wedge (1-a)) \log^{1-1/p} \frac{1}{(a \wedge (1-a))}.$$

We conclude that

$$\mathfrak{Is}_{\mu_{p,n}}(a) \geq c_3 c_4(a \wedge (1-a)) \log^{1-1/p} \frac{1}{(a \wedge (1-a))}$$

and

$$\mathcal{I}_{\mu_{p,n}}(a) \geq c_3 c_4 a \log^{1-1/p} \frac{1}{a}, \quad 0 < a < \frac{1}{2}.$$

*Remark 3.* The proof in [3] relies on rather involved semigroup estimates. For $p = 1$, $\mu = \mu_{1,n}$, the inequality (10) was proved earlier by Bobkov and Houdré [11], using a more elementary argument. For $p = 2$, $\mu = \mu_{2,n}$, (10) follows from the Gaussian isoperimetric inequality that was proved by Sudakov and Tsirelson [26] and Borell [14]; an elementary proof was given by Bobkov [8].

If the Lipschitz semi-norm of $T$ were of order $n^{-1/p}$, the main inequality (2) would follow immediately (since Lipschitz maps preserve isoperimetric inequalities). Unfortunately, $\|T\|_{\text{Lip}} = +\infty$ (as follows from the computation in the proof of Lemma 1). Therefore we use a cut-off argument, cutting off the parts of the space where the local Lipschitz norm is too big. This appears more natural when the isoperimetric inequality is written in a functional form.

Functional forms of isoperimetric inequalities were introduced around 1960 by Maz'ya [22], Federer and Fleming [17]. We follow the approach developed by Bobkov and Houdré [11]; the reader may refer to the latter work for a more general and detailed exposition.

**Proposition A (Maz'ya, Federer–Fleming, Bobkov–Houdré).** *Let $\mu$ be a probability measure, $0 < a < 1/2$, $b > 0$. The following are equivalent:*

(1) $\mathcal{I}_\mu(a) \geq b$;
(2) $\mathcal{J}_\mu(a) \geq b$;
(3) *for any locally Lipschitz function $\phi : \operatorname{supp}\mu \to [0,1]$ such that $\mu\{\phi = 0\} \geq 1/2$ and $\mu\{\phi = 1\} \geq a$,*

$$\int \|\nabla \phi\|_2 \, d\mu \geq b.$$

**Proof.** We shall only prove (1) $\Leftrightarrow$ (3); the proof is similar for (2) $\Longleftrightarrow$ (3).

(1) $\Rightarrow$ (3): by the co-area inequality (which is just the first inequality in the following, cf. Bobkov and Houdré [11]),

$$\int \|\nabla \phi\|_2 \, d\mu \geq \int_0^1 \mu^+\{\phi > u\} \, du \geq \int_0^1 \mathcal{I}_\mu(a) \, du = \mathcal{I}_\mu(a) \geq b.$$

(3) $\Rightarrow$ (1): for a set $A$ of measure $a \leq \mu(A) < 1/2$, let

$$\phi(x) = 0 \vee (1 - s^{-1} \operatorname{dist}(x, \{y \mid \operatorname{dist}(y, A) \leq r\})).$$

If $x \in A$, $\phi(x) = 1$; thereby $\mu\{\phi = 1\} \geq a$. If $\operatorname{dist}(x, A) > r + s$, then $\phi(x) = 0$; the set of all these $x$ has measure $\geq 1/2$ for sufficiently small $r + s$, except maybe for the trivial case $\mu^+(A) = +\infty$.

Further, $\|\nabla \phi(x)\|_2 \leq s^{-1}$, and is 0 unless $r \leq \operatorname{dist}(x, A) \leq r + s$. Therefore according to the assumption

$$s^{-1} \mu\{r \leq \operatorname{dist}(x, A) \leq r + s\} \geq b$$

(for sufficiently small $r + s$). Letting $r, s \to +0$ so that $r/s \to 0$ we recover (1). □

Now let us formulate (and prove) some technical lemmata. First, we need an estimate on the operator norm of the (adjoint) derivative

$$D^* T(z) : \mathbb{R}^n \to \mathbb{R}^{n+1}.$$

**Lemma 1.** *For $0 \neq z \in \mathbb{R}^n$,*

$$\|D^* T(z)\| \leq \frac{1}{\|z\|_p} \{1 + n^{(2-p)/(2p)} \|T(z)\|_2\}.$$

**Proof.** To simplify the notation, assume that $z_i \geq 0$, $i = 1, 2, \ldots, n+1$. Then

$$\frac{\partial T_j}{\partial z_i}(z) = \frac{1}{\|z\|_p} \left\{ \delta_{ij} - \frac{x_j z_i^{p-1}}{\|z\|_p^p} \right\} = \frac{1}{\|z\|_p} \left( \mathbf{1} - \frac{z^{p-1} \otimes x}{\|z\|_p^p} \right)_{ij}$$

(where the power $z^{p-1}$ is computed coordinate-wise). Now,

$$\|z^{p-1} \otimes x\| = \|z^{p-1}\|_2 \times \|x\|_2 = \|z\|_{2(p-1)}^{p-1} \times \|x\|_2$$

$$\leq n^{(2-p)/(2p)} \times \|z\|_p^{p-1} \times \|x\|_2 = n^{(2-p)/(2p)} \times \|z\|_p^p \times \|T(z)\|_2$$

by the Hölder inequality.   □

The following trivial lemma justifies the cut-off arguments.

**Lemma 2.** *If $k, h : \mathbb{R}^n \to [0, 1]$ are two locally Lipschitz functions, then*

$$\|\nabla k\|_2 \geq \|\nabla(kh)\|_2 - \|\nabla h\|_2.$$

Define two cut-off functions

$$\begin{cases} h_1 : \mathbb{R}^n \to [0, 1], & x \mapsto 0 \vee (1 \wedge (2 - c_1 n^{(2-p)/(2p)} \|x\|_2)) \quad \text{and} \\ h_2 : \mathbb{R}^{n+1} \to [0, 1], & z \mapsto 0 \vee (1 \wedge (c_2 n^{-1/p} \|z\|_p - 1)). \end{cases}$$

The function $h_1$ will be used to cut off the part of the space where $\|Tz\|_2$ is too large; the function $h_2$ will be used to cut off the part of the space where $\|z\|_p$ is too small. We shall choose $c_1$ and $c_2$ later on, in the proof of Proposition 2.

The next lemma collects the properties of $h_1$ and $h_2$.

**Lemma 3.** *The function $h_1$ is identically 0 on $\{\|x\|_2 \geq 2c_1^{-1} n^{-(2-p)/(2p)}\}$ and 1 on $\{\|x\|_2 \leq c_1^{-1} n^{-(2-p)/(2p)}\}$. The gradient modulus $\|\nabla h_1\|_2$ is not greater than $c_1 n^{(2-p)/p}$, and vanishes outside*

$$\{c_1^{-1} n^{-(2-p)/(2p)} \leq \|x\|_2 \leq 2c_1^{-1} n^{-(2-p)/(2p)}\}.$$

*The function $h_2$ is 0 on $\{\|z\|_p \leq c_2^{-1} n^{1/p}\}$ and 1 on $\{\|z\|_p \geq 2c_2^{-1} n^{1/p}\}$; $\|\nabla h_2\|_2$ is not greater than $c_2 n^{-1/2}$, and vanishes outside*

$$\{c_2^{-1} n^{1/p} \leq \|z\|_p \leq 2c_2^{-1} n^{1/p}\}.$$

**Proof.** The inequality $\|\nabla h_2\|_2 \leq c_2 n^{-1/2}$ follows from Hölder's inequality:

$$\|h_2(z) - h_2(z')\|_2 \leq c_2 n^{-1/p} \|z - z'\|_p \leq c_2 n^{-1/2} \|z - z'\|_p;$$

the other statements are obvious.   □

Finally, we have the following lemma.

**Lemma 4.** *For any $C_1 > 0$ there exists $C_2 > 0$ (independent of $p \in [1, 2]$ and $n \in \mathbb{N}$) such that*

$$V_{p,n}\{\|x\|_2 \geq C_2 n^{-(2-p)/(2p)}\} \leq \exp(-C_1 n^{p/2})$$

*and*

$$\mu_{p,n}\{\|z\|_p \leq C_2^{-1} n^{1/p}\} \leq \exp(-C_1 n) \leq \exp(-C_1 n^{p/2}).$$

**Proof.** The first part follows from the Schechtman–Zinn theorem (7).

As for the second part,

$$\mu_{p,n}\{\|z\|_p \leq (cn)^{1/p}\} = \mu_{p,n}\Big\{\sum |z_i|^p \leq cn\Big\}. \tag{12}$$

If $Z = (Z_1, \ldots, Z_{n+1}) \sim \mu_{p,n}$, then $|Z_i|^p$ are nonnegative independent random variables. The density of $Z_i$ ($1 \leq i \leq n$) is

$$\frac{x^{-(p-1)/p} \exp(-x)}{\Gamma(1/p)} \mathbf{1}_{[0,+\infty)}(x) \, dx,$$

the density of $Z_{n+1}$ is

$$\exp(-x)\mathbf{1}_{[0,+\infty)}(x)\,\mathrm{d}x,$$

and both are bounded by $\mathrm{const}\cdot x^{-(p-1)/p}\,\mathrm{d}x$ (what is essential here is that the density does not grow too fast near 0). Thus, an estimate on (12) follows from standard large deviation arguments that we reproduce for completeness in Lemma 5. □

**Lemma 5.** *Let $X_1,\ldots,X_N \geq 0$ be independent random variables such that the density of every one of them is bounded by $Ax^{-\alpha}\,\mathrm{d}x$ for some $A>0$ and $0\leq\alpha<1$. Then*

$$\mathbb{P}\{X_1+\cdots+X_N \leq N\varepsilon\} \leq [C(A,\alpha)\varepsilon]^{(1-\alpha)N},$$

*where $C(A,\alpha) = \frac{\mathrm{e}}{1-\alpha}[A\Gamma(1-\alpha)]^{1/(1-\alpha)}$.*

**Proof.** Let $Y = X_1+\cdots+X_N$. For $1\leq i\leq N$, $\xi\geq 0$,

$$\mathbb{E}\exp(-\xi X_i) \leq A\int_0^\infty \exp(-\xi x)x^{-\alpha}\,\mathrm{d}x = A\Gamma(1-\alpha)\xi^{-(1-\alpha)};$$

therefore

$$\mathbb{E}\exp(-\xi Y) \leq [A\Gamma(1-\alpha)]^N \xi^{-N(1-\alpha)}.$$

By Chebyshev's inequality

$$\begin{aligned}\mathbb{P}\{Y\leq N\varepsilon\} &= \mathbb{P}\left\{\exp\left(-\frac{1-\alpha}{\varepsilon}Y\right)\geq \exp(-(1-\alpha)N)\right\}\\ &\leq \exp((1-\alpha)N)\mathbb{E}\exp\left(-\frac{1-\alpha}{\varepsilon}Y\right)\\ &\leq \exp((1-\alpha)N)[A\Gamma(1-\alpha)]^N \left(\frac{1-\alpha}{\varepsilon}\right)^{-(1-\alpha)N}.\end{aligned}$$ □

**Proof of Proposition 2.** Let $0<a<1/2$. Pick $f\colon B_p^n \to [0,1]$ such that $V_{p,n}\{f=0\}\geq 1/2$ and $V_{p,n}\{f=1\}\geq a\geq \exp(-Cn^{p/2})$. Then (by Lemmata 2 and 3)

$$\begin{aligned}\int\|\nabla f\|_2\,\mathrm{d}V_{p,n} &\geq \int\|\nabla(fh_1)\|_2\,\mathrm{d}V_{p,n} - \int\|\nabla h_1\|_2\,\mathrm{d}V_{p,n}\\ &\geq \int\|\nabla(fh_1)\|_2\,\mathrm{d}V_{p,n} - c_1 n^{(2-p)/(2p)} V_{p,n}\{\|x\|_2 \geq c_1^{-1}n^{-(2-p)/(2p)}\}.\end{aligned} \tag{13}$$

Let $g = (fh_1)\circ T$. By the definition of push-forward and Lemma 1,

$$\begin{aligned}\int\|\nabla(fh_1)\|_2\,\mathrm{d}V_{p,n} &= \int\|\nabla(fh_1)\circ T\|_2\,\mathrm{d}\mu_{p,n}\\ &\geq \int_{\mathbb{R}^{n+1}} \frac{\|\nabla g(z)\|_2}{\|D^*T(z)\|}\,\mathrm{d}\mu_{p,n}\\ &\geq \int_{\mathbb{R}^{n+1}} \frac{\|\nabla g(z)\|_2\|z\|_p}{1+n^{(2-p)/(2p)}\|T(z)\|_2}\,\mathrm{d}\mu_{p,n}.\end{aligned}$$

According to Lemma 3, $\|T(z)\|_2 \leq 2c_1^{-1}n^{-(2-p)/(2p)}$ whenever $h_1(T(z))\neq 0$; hence

$$\int\|\nabla(fh_1)\|_2\,\mathrm{d}V_{p,n} \geq c_3\int_{\mathbb{R}^{n+1}}\|\nabla g\|_2\|z\|_p\,\mathrm{d}\mu_{p,n}, \tag{14}$$

where $c_3 = c_1/(c_1 + 2)$. Applying Lemmata 2 and 3 once again,

$$\int_{\mathbb{R}^{n+1}} \|\nabla g\|_2 \|z\|_p \, d\mu_{p,n} \geq \int_{\mathbb{R}^{n+1}} \|\nabla(gh_2)\|_2 \|z\|_p \, d\mu_{p,n} - \int_{\mathbb{R}^{n+1}} \|\nabla h_2\|_2 \|z\|_p \, d\mu_{p,n}$$

$$\geq c_2^{-1} n^{1/p} \int_{\mathbb{R}^{n+1}} \|\nabla(gh_2)\|_2 \, d\mu_{p,n} - 2n^{(2-p)/(2p)} \mu_{p,n}\{\|z\|_p \leq 2c_2^{-1} n^{1/p}\}. \quad (15)$$

The inequalities (13)–(15) show that

$$\int \|\nabla f\|_2 \, dV_{p,n} \geq c_4 n^{1/p} \int_{\mathbb{R}^{n+1}} \|\nabla(gh_2)\|_2 \, d\mu_{p,n} - c_1 n^{(2-p)/(2p)} V_{p,n}\{\|x\|_2 \geq c_1^{-1} n^{-(2-p)/(2p)}\}$$

$$- 2c_3 n^{(2-p)/(2p)} \mu_{p,n}\{\|z\|_p \leq 2c_2^{-1} n^{1/p}\}, \quad (16)$$

where $c_4 = c_3 c_2^{-1}$. Therefore by Lemma 4 (with $C_1$ larger than $C$ from Proposition 1) we can choose $c_1$ and $c_2$ so that

$$\int \|\nabla f\|_2 \, dV_{p,n} \geq c_4 n^{1/p} \int_{\mathbb{R}^{n+1}} \|\nabla(gh_2)\|_2 \, d\mu_{p,n} - \exp(-Cn^{p/2})/2. \quad (17)$$

The function $gh_2 = ((fh_1) \circ T)h_2$ vanishes on a set of $\mu_{p,n}$-measure $\geq 1/2$ (for example, it is zero on $T^{-1}\{f = 0\}$). Also,

$$\{gh_2 = 1\} \supset \{g = 1\} \setminus \{h_2 < 1\} \supset T^{-1}(\{f = 1\} \setminus \{h_1 < 1\}) \setminus \{h_2 < 1\}$$

is of $\mu_{p,n}$-measure at least

$$V_{p,n}\{f = 1\} - V_{p,n}\{\|x\|_2 > c_1^{-1} n^{-(2-p)/(2p)}\} - \mu_{p,n}\{\|z\|_2 < 2c_2^{-1} n^{1/p}\}$$

$$\geq V_{p,n}\{f = 1\} - \exp(-Cn^{p/2})/2 \geq \frac{1}{2} V_{p,n}\{f = 1\} \geq \frac{a}{2}.$$

Therefore by inequality (10) and Proposition A

$$\int_{\mathbb{R}^{n+1}} \|\nabla(gh_2)\|_2 \, d\mu_{p,n} \geq c \frac{a}{2} \log^{1-1/p} \frac{2}{a} \geq c_5 a \log^{1-1/p} \frac{1}{a}.$$

To conclude, combine this inequality with (17) and apply Proposition A once again. □

## 4. Remarks

(1) Let us briefly recall the connection between the isoperimetric inequality as in Theorem 1 and related $L_2$ inequalities.

According to Proposition A, (2) for $\mu = V_{p,n}$ can be written as

$$\int \|\nabla \phi\|_2 \, d\mu \geq cn^{1/p} a \log^{1-1/p} \frac{1}{a} \quad \text{for } 0 \leq \phi \leq 1 \text{ such that } \mu\{\phi = 0\} \geq \frac{1}{2}, \ \mu\{\phi = 1\} \geq a. \quad (18)$$

The following is well known.

**Proposition B.** *If a probability measure $\mu$ satisfies (18), then also*

$$\int \|\nabla \phi\|_2^2 \, d\mu \geq c_1 n^{2/p} a \log^{2-2/p} \frac{1}{a} \quad \text{for } 0 \leq \phi \leq 1 \text{ such that } \mu\{\phi = 0\} \geq \frac{1}{2}, \ \mu\{\phi = 1\} \geq a \quad (19)$$

*(with some constant $c_1$ depending on $c$).*

As proved by Barthe and Roberto [5], (19) is (up to constants and normalisation) an equivalent form of the Latała–Oleszkiewicz inequality (introduced in [21] under the name $I(p)$).

**Proof of Proposition B.** Assume for simplicity that $\phi$ has no atoms except for 0 and 1. For $0 \le u, \varepsilon \le 1$, let $\phi_{u,\varepsilon} = 0 \vee (1 \wedge \varepsilon^{-1}(\phi - u))$. By (18) and Jensen's inequality,

$$\int_{u \le \phi \le u+\varepsilon} \|\nabla \phi\|_2^2 \, d\mu = \varepsilon^2 \int \|\nabla \phi_{u,\varepsilon}\|_2^2 \, d\mu$$
$$\ge \frac{\varepsilon^2}{\mu\{u \le \phi \le u+\varepsilon\}} \left[ \int \|\nabla \phi_{u,\varepsilon}\|_2 \, d\mu \right]^2$$
$$\ge \frac{\varepsilon^2}{\mu\{u \le \phi \le u+\varepsilon\}} c^2 n^{2/p} m_{u+\varepsilon}^2 \log^{2-2/p} \frac{1}{m_{u+\varepsilon}},$$

where $m_u = \mu\{\phi \ge u\}$. Let $u_0 = 0$, $u_{i+1} = u_i + \varepsilon_i$, choosing $\varepsilon_i$ so that

$$\mu\{u_i < \phi \le u_i + \varepsilon_i\} = m_{u_i + \varepsilon_i} = 1/2^{i+1}$$

(except for the last step, $i = [\log 1/a]$). As $\sum \varepsilon_i = 1$, the Cauchy–Schwarz inequality yields

$$\int \|\nabla \phi\|_2^2 \, d\mu \ge \sum_i c^2 \varepsilon_i^2 n^{2/p} m_{u_{i+1}} \log^{2-2/p} \frac{1}{m_{u_{i+1}}} \ge c^2 n^{2/p} \left( \sum_{1 \le i \le \log(1/a)} 2^i \log^{-(2-2/p)} 2^i \right)^{-1}$$
$$\ge c_1 n^{2/p} a \log^{2-2/p} \frac{1}{a}. \qquad \square$$

In the class of log-concave measures, the last proposition can be reversed (that is, (19) implies (18) with a constant $c$ depending on $c_1$). This was proved by Michel Ledoux [20] for $p = 1, 2$ (see also [6], Theorem 1.3), and extended by Barthe–Cattiaux–Roberto [3] to all $1 \le p \le 2$.

(2) The volume measure $V_K = \text{Vol}|_K / \text{Vol}(K)$ on a convex body $K \subset \mathbb{R}^n$ has attracted much interest in recent years. For any body $K \subset \mathbb{R}^n$ there exists a nondegenerate linear map $T : \mathbb{R}^n \to \mathbb{R}^n$ such that $\widetilde{K} = TK$ is *isotropic*, meaning that

$$\text{Vol}\,\widetilde{K} = 1, \qquad \int_{\widetilde{K}} x_i x_j \prod_{k=1}^n dx_k = L_K^2 \delta_{ij} \quad \text{for } 1 \le i, j \le n.$$

The number $L_K$ is an invariant of the body called the isotropic constant; $L_K > c$ for some universal constant $c > 0$ (independent of $K$ and $n$). The famous hyperplane conjecture asserts that $L_K \le C$. So far, it is only known that $L_K \le Cn^{1/4}$; this was recently proved by Bo'az Klartag [18], improving the celebrated estimate of Bourgain [15] with an extra logarithmic factor.

Kannan, Lovász and Simonovits [19] conjectured that there exists a universal constant $c_0 > 0$ such that for any isotropic convex body $K$ the measure $\mu = V_K$ satisfies the Cheeger-type inequality

$$\mathcal{I}_\mu(a) \ge \frac{c_0 a}{L_K} \tag{20}$$

for $0 < a \le 1/2$.

The inequality (20) has so far been proved for a mere few families of convex bodies (cf. [9, 12] for an extensive discussion and related results, and [7, 11] for several families of examples in the larger class of log-concave measures). As $\widetilde{B_p^n} = C(n,p) B_p^n$, where $cn^{1/p} \le C(n,p) \le Cn^{1/p}$, Theorem 1 shows that the conjecture is true for $\widetilde{B_p^n}$, $1 \le p \le 2$.

Recently, Grigoris Paouris [23] proved that for any isotropic convex body $K$,

$$V_K\{\|x\|_2 \ge t\} \le \exp\left(-\frac{ct}{L_k}\right), \quad t \ge t_0 L_K \sqrt{n}. \tag{21}$$

Repeating the proof of Proposition 1 with (21) instead of the Schechtman–Zinn theorem, we obtain:

**Proposition B.** *If $K \subset \mathbb{R}^n$ is an isotropic convex body, then* (20) *holds for $0 < a < \exp(-C\sqrt{n})$ (where $C$ is a universal constant).*

(3) The right-hand side in the inequality (2) behaves like $c(n) a \log^{1/p} \frac{1}{a}$ as $a \to 0$, whereas the correct asymptotics should be $c(n) a^{1-1/n}$ (the difference becomes essential however only for $a \lesssim e^{-cn \log n}$).

To recover the correct behaviour for small $a$, note that the inequality (9) that we used in the proof of Proposition 1 is dimension free. We can use instead the following dimensional extension, due to Franck Barthe [2]:

**Theorem (Barthe).** *Let $K$ be a convex body in $\mathbb{R}^n$ and let $V_K$ be the normalised volume measure on $K$. Then, for any $A \subset K$ and any $r > 0$,*

$$V_K^+(A) \geq \frac{n}{2r} \{[V_K(A)^{1-1/n} + (1 - V_K(A))^{1-1/n}] V_K\{\|x\|_2 \leq r\}^{1/n} - 1\}.$$

## Acknowledgments

I thank Franck Barthe for illuminating discussions, for sharing his interest in functional inequalities, and for bringing the problem to my attention. I thank my supervisor Vitali Milman for encouragement and support. I thank Sergey Bobkov for comments and explanations.

I thank them and the anonymous referees for the remarks on preliminary versions of this note.

This work was done while the author enjoyed the hospitality of the Paul Sabatier University, Toulouse, staying there on a predoc position of the European network PHD.